\documentclass{article}
\usepackage{amssymb,amsmath,mathrsfs,latexsym,amscd}
\usepackage{exscale,latexsym,amsthm,graphics}

\usepackage{makeidx}
\makeindex
\usepackage{epsfig}
\usepackage{txfonts}
\newtheorem{thm}{Theorem}[section]
\newtheorem{cor}[thm]{Corollary}

\newtheorem{lemma}[thm]{Lemma}

\newtheorem{proposition}[thm]{Proposition}

\newtheorem{remark}[thm]{Remark}

\def\qed{{\hfill $\Box$ \bigskip}}

\def\R {{\mathbb R}}

\def\P{{\mathbb P}}

\def\E{{\mathcal E}}
\def\P{{\mathbb P}}

\numberwithin{equation}{section}
\begin{document}

\noindent
{{\Large\bf  Construction of Liouville Brownian motion via Dirichlet form theory}}

\bigskip
\noindent
{\bf Jiyong Shin}
\\

\noindent
{\small{\bf Abstract.} 
The Liouville Brownian motion which was introduced in \cite{GRV} is a natural diffusion process associated with a random metric in two dimensional Liouville quantum gravity. In this paper we construct the Liouville Brownian motion via Dirichlet form theory. By showing that the Liouville measure is smooth in the strict sense, the positive continuous additive functional $(F_t)_{t \ge 0}$ of the Liouville measure in the strict sense w.r.t.  the planar Brownian motion $(B_t)_{t \ge 0}$ is obtained. Then the Liouville Brownian motion can be defined as a time changed process of the planar Brownian motion $B_{F_t^{-1}}$.\\



\section{Introduction}

In \cite{GRV} C. Garban, R. Rhodes, V. Vargas constructed the Liouville Brownian motion (hereafter LBM) which is a canonical diffusion process in planar Liouville quantum gravity. By classical theory of Gaussian multiplicative chaos (cf. \cite{Kah}), the Liouville measure $M_{\gamma}$, $\gamma \in (0,2)$ (see Section \ref{ss;ml}) is formally defined as 
\[
M_{\gamma} (dz)=  \exp \Big( \gamma X(z) - \frac{\gamma^2}{2} E[X(z)^2] \Big) \ dz,
\]
where $X$ is a massive Gaussian free field on $\R^2$ on a probability space $(\Omega,\mathcal{A},P)$ and $dz$ is the Lebesgue measure on $\R^2$. By  constructing the positive continuous additive functional $(A_t)_{t \ge 0}$ of a Brownian motion $(W_t)_{t \ge 0}$ w.r.t. the Liouville measure $M_{\gamma}$, the Liouville Brownian motion is defined as $W_{A_t^{-1}}$.  Subsequently, the heat kernel estimates and Dirichlet forms associated with the LBM were investigated in \cite{GRV2} and \cite{AnKa}. 

In this short paper, we are concerned with constructing the Liouville Brownian motion via Dirichlet form theory.  It is known from \cite{GRV2} and \cite{FOT} that the LBM is associated with the Liouville Dirichlet form. Therefore it is a natural question if one can construct the Liouville  Brownian motion in view of Dirichlet form theory \cite{FOT} (cf. \cite[Section 1.4]{GRV2}). In general there is no theory of  Dirichlet forms which enables to get rid of the polar set and construct a Hunt process starting from all points of $\R^2$. Recently, using elliptic regularity results only n-regularized Liouville Brownian motion was constructed via Dirichlet form theory (see \cite{Sh}). This approach, however, can not be applied to construct the LBM since the massive Gaussian free field is not a function. In principle our construction of LBM is based on the observation of power law  decay of the size of balls of the Liouville measure and Proposition \ref{p;smoothl}. By Lemma \ref{l;contigr} and Proposition \ref{p;smoothl} we can show that $M_{\gamma} \in S_1$ (see Section \ref{s;dfth} for the definition of $S_1$).  Therefore there exists a positive continuous additive functional $(F_t)_{t \ge 0}$ in the strict sense of the planar Brownian motion $(B_t)_{t \ge 0}$ w.r.t. $M_{\gamma}$ by Revuz correspondence (see \cite[Theorem 5.1.7]{FOT}). Then the Liouville Brownian motion can be defined as $B_{F_t^{-1}}$.\\

\noindent{\bf{Notations}}:\\
We denote the set of all  Borel measurable functions and the set of all bounded Borel measurable functions on $\R^2$  by  $\mathcal{B}(\R^2)$ and $\mathcal{B}_b(\R^2)$, respectively. The usual $L^q$-spaces $L^q(\R^2, dx)$, $q \in[1,\infty]$ are equipped with $L^{q}$-norm $\| \cdot \|_{L^q (\R^2,dx) }$ with respect to the  Lebesgue measure $dx$ on $\R^2$ and $\mathcal{A}_b$ : = $\mathcal{A} \cap \mathcal{B}_b(\R^2)$ for $\mathcal{A} \subset L^q(\R^2,dx)$.  The inner product on $L^2(\R^2, dx)$ is denoted by $(\cdot,\cdot)_{L^2(\R^2, dx)}$.  The indicator function of a set $A$ is denoted by $1_A$. Let $\nabla f : = ( \partial_{1} f, \dots , \partial_{d} f )$  and  $\Delta f : = \sum_{j=1}^{d} \partial_{jj} f$ where $\partial_j f$ is the $j$-th weak partial derivative of $f$ and $\partial_{jj} f := \partial_{j}(\partial_{j} f) $, $j=1, \dots, d$. The Sobolev space $H^{1,q}(\R^2, dx)$, $q \ge 1$ is defined to be the set of all functions $f \in L^{q}(\R^2, dx)$ such that $\partial_{j} f \in L^{q}(\R^2, dx)$, $j=1, \dots, d$, and $H^{1,q}_{loc}(\R^2, dx) : =  \{ f  \,|\;  f \cdot \varphi \in H^{1,q}(U, dx),\,\forall \varphi \in  C_0^{\infty}(U)\}$. 
Here $C_0^{\infty}(\R^2)$ denotes the set of all infinitely differentiable functions with compact support in $\R^2$. We also denote the set of continuous functions on $\R^2$ and  the set of compactly supported continuous functions on $\R^2$ by $C(\R^2)$ and $C_0(\R^2)$, respectively. We equip $\R^2$ with the Euclidean norm $|\cdot|$ and the corresponding inner product $\langle \cdot, \cdot \rangle$.

\section{Preliminary}
\subsection{Dirichlet form theory}\label{s;dfth}
In this subsection we present some definitions and properties in Dirichlet form theory as stated in \cite{FOT}. Let 
\[
\E(f,g):= \frac{1}{2} \int_{\R^2}   \langle \nabla f,  \nabla g \rangle \ dx, 
\]
where $f, g \in D(\E) = \{f \in L^2(\R^2,dx) \mid \nabla f \in L^2(\R^2,dx) \}$. Let  $( (B_t)_{t \ge 0}, (\P_x)_{x \in \R^2})$ be the planar Brownian motion associated to the Dirichlet form $(\E,D(\E))$. A positive Radon measure $\mu$ on $\R^2$ is said to be of finite energy integral if
\[
\int_{\R^2} |f(x)|\, \mu (dx) \leq c \sqrt{\E_1(f,f)}, \; f\in D(\E) \cap C_0(\R^2),
\]
where $c$ is some constant independent of $f$ and $\E_{1}(\cdot,\cdot):=\E(\cdot,\cdot) + (\cdot,\cdot)_{L^2(\R^2,dx)}$.  A positive Radon measure $\mu$ on $\R^2$ is of finite energy integral if and only if there exists a unique function $U_{1} \ \mu\in D(\E )$ such that
\[
\E_{1}(U_{1} \, \mu, f) = \int_{\R^2} f(x) \, \mu(dx),
\]
for all $f \in D(\E) \cap C_0(\R^2)$. $U_{1} \, \mu$ is called $1$-potential of $\mu$. For any set $A \subset \R^2$ the capacity of $A$ is defined as
\[
\text{Cap(A)} = \inf_{\begin{subarray} \ \text{open} \ B \subset \R^2  \\    A \subset B \end{subarray}} \inf_{\begin{subarray} \  f \in D(\E) \\ 1_{B} \cdot f \ge 1 \ \text{dx-a.e.}  \end{subarray}} \E_1(f,f).
\]
We denote by $p_{t}(x,y)$ the transition kernel density of $(B_t)_{t \ge 0}$. Taking the Laplace transform of $p_{\cdot}(x, y)$, we see that there exists a $\mathcal{B}(\R^2) \times \mathcal{B}(\R^2)$ measurable non-negative map $r_{1}(x,y)$ such that
\[
R_1 f(x) := \int_{\R^2} r_1 (x,y) \ f(y) \ dy \quad  x \in \R^2, \ f \in \mathcal{B}_b(\R^2).
\]
The density $r_1 (x,y)$ is called the resolvent kernel density. For a signed Radon measure $\mu$ on $\R^2$, let us define
\[
R_1 \mu (x) = \int_{\R^2} r_1 (x,y) \ \mu(dy), \quad  x \in \R^2,
\]
whenever this makes sense.  In particular, $R_{1} \mu$ is a version of $U_{1} \mu$ (see e.g. \cite[Exercise 4.2.2]{FOT}). The family of the measures of finite energy integral is denoted by $S_0$. We further define $S_{00} : = \{\mu\in S_0 \mid \mu(\R^2)<\infty, \|U_{1} \mu\|_{L^{\infty}(\R^2,dx)}<\infty \}$.
A positive Borel measure $\mu$ on $\R^2$ is said to be smooth in the strict sense if there exists a sequence $(E_k)_{k\ge 1}$ of Borel sets increasing to $\R^2$ such that $1_{E_{k}} \cdot \mu \in S_{00}$ for each $k$ and 
\[
\P_{x} ( \lim_{k \rightarrow \infty} \sigma_{ \R^2 \setminus E_{k} }  \ge \zeta ) =1, \quad \forall x \in \R^2,
\]
where $\zeta$ is the lifetime of $(B_t)_{t \ge 0}$.
The totality of the smooth measures in the strict sense is denoted by $S_{1}$ (see \cite{FOT}). If  $\mu \in S_{1}$,
then there exists a unique $A \in A_{c,1}^{+}$ with $\mu = \mu_{A}$, i.e. $\mu$ is the Revuz measure of $A$ (see \cite[Theorem 5.1.7]{FOT}), where $A_{c,1}^{+}$ denotes the family of the positive continuous additive functionals on $\R^2$ in the strict sense.
\subsection{Massive Gaussian free field and Liouville measure}\label{ss;ml}
The massive Gaussian free field on $\R^2$ is a centered Gaussian random distribution (in the sense of Schwartz)  on a probability space $(\Omega, \mathcal{A}, P)$ with covariance function given by the Green function $G^{(m)}$ of the operator $m^2 -\Delta$, $m > 0$, i.e.
\[
(m^2 - \Delta) G^{(m)}(x, \cdot) = 2 \pi \delta_x, \quad x \in \R^2, 
\] 
where $\delta_x$ stands for the Dirac mass at $x$.
The massive Green function with the operator $(m^2 - \Delta)$ can be written as 
\[
G^{(m)}(x,y) = \int_0^{\infty} e^{- \frac{m^2}{2}s - \frac{\| x-y \|^2}{2s}} \frac{ds}{2s} 
= \int_{1}^{\infty} \frac{k_{m} (s (x-y)) }{s} \ ds, \quad x,y \in \R^2,
\]
where 
\[
k_m (z) = \frac{1}{2} \int_0^{\infty} e^{- \frac{m^2}{2s} \| z  \|^2 - \frac{s}{2} } \ ds.
\]

Note that this massive Green function is  a kernel of $\sigma$-positive type in the sense of Kahane \cite{Kah} since we integrate a continuous function of positive type w.r.t. a positive measure.
Let $(c_n)_{n \ge 1}$ be an unbounded strictly increasing sequence such that $c_1 = 1$ and  $(Y_n)_{n \ge 1}$ be a family of  independent centered continuous Gaussian fields on $\R^2$ on the probability space $(\Omega, \mathcal{A}, P)$ with covariance kernel given by 
\[
E[Y_n (x) \ Y_n (y)] = \int_{c_{n-1}}^{c_n}  \frac{k_m (s(x-y))}{s} ds.
\]

The massive Gaussian free field is the Gaussian distribution defined by 
\[
X(x) = \sum_{k \ge 1} Y_k (x).
\]
We define $n$-regularized field by
\[
X_n(x) = \sum_{k=1}^{n} Y_k(x), \quad n \ge 1
\]
and the associated  $n$-regularized Liouville measure by
\[
M_{n,\gamma} (dz)=  \exp \Big( \gamma X_n(z) - \frac{\gamma^2}{2} E[X_n(z)^2] \Big) \ dz, \quad \gamma \in (0,2).
\]
By the classical theory of Gaussian multiplicative chaos (see \cite{Kah}), $P$-a.s. the family $(M_{n,\gamma})_{n \ge 1}$ weakly converges to the measure $M_{\gamma}$, which is called Liouville measure. It is known from \cite{Kah} that $M_{\gamma}$ is a Radon measure on $\R^2$ and has full support.

\section{Liouville Brownian motion via Dirichlet form theory}
In this section we show that the LBM can be constructed in the framework of \cite{FOT}. The following lemma is a direct consequence of \cite[Proposition 2.3]{GRV}: 
\begin{lemma}\label{l;contigr}
Almost surely in $X$, the mapping
\[
x \longmapsto \int_{\R^2} \log_{+}  \frac{1}{ |x -y |} \ M_{\gamma}(dy)
\]
is continuous on $\R^2$ where $\log_{+} a : = \max \{\log a , 0 \}$.
\end{lemma}
The estimates of resolvent kernel density of the planar Brownian motion $(B_t)_{t \ge 0}$ are well known:
\begin{lemma}\label{l;rdel} 
For any $x,y \in \R^2$
\[
r_1(x,y) \le \  c_1 \log_{+} \frac{1}{| x-y |} + c_2, 
\]
where $c_1, c_2>0$ are some constants.
\end{lemma}
\begin{lemma}\label{l;s0}
Almost surely in $X$, for any relatively compact open set $G$, $1_{G} \cdot M_{\gamma} \in S_{0}$.
\end{lemma}
\proof
Note that by Lemma \ref{l;contigr} and Lemma \ref{l;rdel}
\[
\int_G \int_G r_1(x,y) \ M_{\gamma}(dy) \ M_{\gamma}(dx) < \infty. 
\]
Hence $1_G \cdot M_{\gamma} \in S_0$ by \cite[Example 4.2.2]{FOT}
\qed

Now we restate  \cite[Proposition 2.13]{ShTr13a} in our setting, which is a main ingredient to construct the LBM:
\begin{proposition}\label{p;smoothl}
Let $\mu$ be a positive Radon measure on $\R^2$. Suppose that for some relatively compact open set $G \subset \R^2$, $1_{G} \cdot \mu \in S_{0}$  and  $R_1(1_{G} \cdot \mu)$ is bounded $dx$-a.e. by a continuous function $r_1^G \in C(\R^2)$. Then $1_{G} \cdot \mu \in S_{00}$. In particular, if this holds for any relatively compact open set $G$, then $\mu \in S_{1}$.
\end{proposition}

\begin{thm}
Almost surely in $X$, the Liouville measure $M_{\gamma} \in S_{1}$.
\end{thm}
\proof
Let $(E_k)_{k\ge 1}$ be an increasing sequence of relatively compact open sets with $\bigcup_{k \ge1} E_k = \R^2$.
Since $M_{\gamma}$ is a Radon measure, $1_{E_k} \cdot M_{\gamma} (\R^2) < \infty$. For any $x \in \R^2$
\begin{eqnarray}\label{eq;calpo1}
R_1 (1_{E_k} \cdot M_{\gamma} ) (x)  &=& \int_{E_k} r_1(x,y) \  M_{\gamma}(dy) \le  \int_{E_k} \left( c_1 \log_{+}  \frac{1}{ |x -y |} +c_2 \right) \ M_{\gamma}(dy) \notag \\
&\le& c_1  \int_{\R^2} \log_{+}  \frac{1}{ |x -y |} \ M_{\gamma}(dy) + c_k,
\end{eqnarray}
where $c_k= c_2 \ M_{\gamma} (E_k)$. By Lemma \ref{l;contigr}, the right hand side of \eqref{eq;calpo1} is continuous on $\R^2$.
Therefore by Proposition \ref{p;smoothl}, $M_{\gamma}  \in S_1$.
\qed

Therefore by \cite[Thorem 5.1.7]{FOT} there exists a positive continuous additive functional $(F_t)_{t \ge 0}$ of the Brownian motion $(B_t)_{t \ge 0}$  in the strict sense w.r.t. $M_{\gamma}$. Finally the LBM can be defined as
\[
\mathcal{B}_t = B_{F_t^{-1}}, \quad t \ge 0,
\] 
where $F_t^{-1}:= \inf \{s>0 \mid F_s  > t  \}$. 
As a by-product of our approach, we can directly obtain that the Liouville measure $M_{\gamma}$ charges no set of zero capacity (see \cite[Lemma 1.5]{GRV2}).
\begin{cor}
Let $ (E_{k})_{k \ge 1}$ be an increasing sequnce of relatively compact open sets with $\bigcup_{k \ge 1} E_k =\R^2$. Suppose $\mu \in S_1$  with respect to the sequence of the sets $ (E_{k})_{k \ge 1}$. Then $\mu$ does not charge capacity zero sets. In particular, $M_{\gamma}$ charges no set of zero capacity. 
\end{cor}
\proof
Note that by \cite[Lemma 2.2.3]{FOT} , $1_{E_k} \cdot  \mu $ charges no set of zero capacity (see  Lemma \ref{l;s0}).
Suppose $N \subset \R^2$ be an open set such that Cap$(N) =0$. Then the statement follows from 
\[ 
\mu(N) \le \sum_{k \ge 1} \mu(1_{E_k} \cap N) = 0.
\]
\qed
\begin{remark}
\begin{itemize}
\item[(i)]   Note that the support and quasi support of $M_{\gamma}$ is the whole space $\R^2$  (see \cite{GRV2}). Since $M_{\gamma}$ charges no set of zero capacity, the LBM is associated to the Liouville Dirichlet form\[
\E(f,g) = \frac{1}{2}  \int_{\R^2} \langle \nabla f, \nabla g \rangle \ dx, \quad f,g \in \mathcal{F},
\]
where $\mathcal{F} = \{ f \in L^2(\R^2, M_{\gamma}) \cap H^{1,2}_{loc} (\R^2, dx) \mid \nabla f \in L^2(\R^2, dx) \}$ (see \cite{FOT}).
\item[(ii)] The LBM at criticality (i.e. $\gamma =2$) was introduced in \cite{RhVa}.
In this case our approach can not be applied to construct the LBM at criticality since  the continuous mapping as in Lemma \ref{l;contigr} for the Liouville measure at criticality can not be obained.  
\end{itemize}
\end{remark}

\vspace*{2cm}
\noindent Jiyong Shin\\
School of Mathematics \\
Korea Institute for Advanced Study\\
85 Hoegiro Dongdaemun-gu,\\
Seoul 02445, South Korea,  \\
E-mail: yonshin2@kias.re.kr \\

\addcontentsline{toc}{chapter}{References}
\begin{thebibliography}{99}

\bibitem{AnKa}
S. Andres, N. Kajino, {\it Continuity and estimates of the Liouville heat kernel with applications to spectral dimensions}, Probab. Theory Related Fields 166 (2016), no. 3-4, 713-752. 

\bibitem{FOT}
M. Fukushima, Y. Oshima, M. Takeda, {\it Dirichlet forms and symmetric Markov processes}, Second revised and extended edition. de Gruyter Studies in Mathematics, 19. Walter de Gruyter Co., Berlin (2011).

\bibitem{GRV}
C. Garban, R. Rhodes, V. Vargas, {\it Liouville Brownian motion}, Ann. Probab. 44 (2016), no. 4, 3076-3110. 
\bibitem{GRV2}
C. Garban, R. Rhodes, V. Vargas, {\it On the heat kernel and the Dirichlet form of Liouville Brownian Motion}, Electron. J. Probab. 19 (2014), no. 96.
\bibitem{Kah}
J. P. Kahane, {\it Sur le chaos multiplicatif}, Ann. Sci. Math. Québec 9 (1985), no. 2, 105-150.
\bibitem{RhVa}
R. Rhodes, V. Vargas, {\it Liouville Brownian motion at criticality}, 
Potential Anal. 43 (2015), no. 2, 149-197. 
\bibitem{Sh}
J. Shin,  {\it Elliptic regularity results: n-regularized Liouville Brownian motion and non-symmetric diffusions associated with degenerate forms}, arXiv:1606.05857.
\bibitem{ShTr13a}
J. Shin, G. Trutnau, {\it On the stochastic regularity of distorted Brownian motions}, Trans. Amer. Math. Soc. 369 (2017), no. 11, 7883-7915.


\end{thebibliography}
\end{document}